\documentclass [a4paper,12pt]{article}
\usepackage{graphicx}
\usepackage[small]{subfigure,epsfig}
\usepackage{indentfirst}
\usepackage{amsmath,latexsym,enumerate}
\usepackage{amssymb}
\usepackage{amsfonts}
\usepackage{xcolor}
\usepackage{authblk}
\usepackage{amsthm}
\usepackage{appendix}
\usepackage{array}
\usepackage{tikz}
\usepackage{cite}

\numberwithin{equation}{section}

\usepackage{geometry} 
\begin{document}

\title{Scenarios for the appearance of strange attractors in a model of three interacting microbubble contrast agents}

\author[1,2]{Ivan Garashchuk}

\author[3]{Alexey Kazakov}

\author[4,5,2]{Dmitry Sinelshchikov}

\affil[1]{BCAM - Basque Center for Applied Mathematics, Mazarredo 14, E48009, Bilbao, Basque Country - Spain}

\affil[2]{HSE University, 34 Tallinskaya Street, 123458, Moscow, Russian Federation}

\affil[3]{HSE University, 25/12 Bolshaya Pecherskaya, 10587, Nizhny Novgorod, 10587,  Russian Federation}

\affil[4]{Instituto Biofisika (UPV/EHU, CSIC), University of the Basque Country, Leioa, 48940, Spain}

\affil[5]{Ikerbasque, Basque Foundation for Science, Bilbao 48009, Spain}

\maketitle

\begin{abstract}
We study nonlinear dynamics in a model of three interacting encapsulated gas bubbles in a liquid. The model is a system of three coupled nonlinear oscillators with an external periodic force. Such bubbles have numerous applications, for instance, they are used as contrast agents in ultrasound visualization. Certain types of bubbles oscillations may be beneficial or undesirable depending on a given application and, hence, the dependence of the regimes of bubbles oscillations on the control parameters is worth studying. We demonstrate that there is a wide variety of types of dynamics in the model by constructing a chart of dynamical regimes in the control parameters space. Here we focus on hyperchaotic attractors characterized by three positive Lyapunov exponents and strange attractors with one or two positive Lyapunov exponents possessing an additional zero Lyapunov exponent, which have not been observed previously in the context of bubbles oscillations. We also believe that we provide a first example of a hyperchaotic attractor with additional zero Lyaponov exponent. Furthermore, the mechanisms of the onset of these types of attractors are still not well studied. We identify two-parametric regions in the control parameter space where these hyperchaotic and chaotic attractors appear and study one-parametric routes leading to them. We associate the appearance of hyperchaotic attractors with three positive Lyapunov exponents with the inclusion of a periodic orbit with a three-dimensional unstable manifold, while the onset of chaotic oscillations with an additional zero Lyapunov exponent is connected to the partial synchronization of bubbles oscillations. We propose several underlying bifurcation mechanisms that explain the emergence of these regimes. We believe that these bifurcation scenarios are universal and can be observed in other systems of coupled oscillators.
\end{abstract}

\section{Introduction}

In this work we study a mathematical model of three interacting gas bubbles encapsulated into viscoelastic shells. Such bubbles of micrometer size are used as contrast agents to enhance ultrasound imaging, where they are known as ultrasound contrast agents, for targeted drug delivery and other applications, e.g. sonoporation \cite{Hoff2001,Goldberg2001,Szabo2004,Klibanov2006,Forbes2012,Lentacker2014,Fournier2023}.

Encapsulated microbubbles can exhibit complicated behavior under the influence of the external ultrasound field. Moreover, the character of the oscillations strongly depends on the parameters and  the initial conditions, while the particular type of dynamics may be either desirable or undesirable depending on an application. For example, it is known that irregular, chaotic oscillations can be beneficial for ultrasound imaging due to high contrast with both responses of the surrounding tissues and incident acoustic field \cite{Hoff2001}. Furthermore, for some applications, such as targeted drug delivery or sonoporation, where a quick controlled rupture or large-amplitude oscillations are desirable, hyperchaotic oscillations also can be beneficial. Consequently, understanding of possible types of bubbles dynamics and their bifurcations under the variation of the control parameters is an important problem from an applied point of view. Moreover, the models of bubbles oscillations can be a source of new universal bifurcation scenarios, as it was shown e.g. in \cite{Garashchuk2019,Garashchuk2021} and, hence, their investigation is interesting from a theoretical point of view. In this work we show that this is also the case for the considered model of three interacting bubbles oscillations.

Mathematical models of the bubbles dynamics are based on the Rayleigh--Plesset equation, which was derived for the description of oscillations of a spherical bubble in an infinite incompressible viscous liquid \cite{Plesset1949}. In order to describe oscillations of ultrasound contrast agents this equation needs generalizations that take into account bubbles shell, liquid's compressibility, interactions between bubbles and other aspects (see, e.g. \cite{Takahira1995,Mettin1997,Ida2002,Allen2003,Doinikov2009,Doinikov2009a,Doinikov2011,Paul2010,Carroll2013,Doinikov2015,Versluis2020} and references therein).

Dynamical systems governing oscillations of isolated or coupled microbubbles have been studied in a number of works, and it has been shown that they exhibit complex dynamics including various bifurcation chains and the corresponding interesting attractors (see, e.g. \cite{Parlitz1990,Macdonald2006,Behnia2009,Carroll2013,Dzaharudin2013,Garashchuk2018,Ohl2019,Garashchuk2019,Garashchuk2020,Garashchuk2021,Sojahrood2021,Sojahrood2021a}). For instance, authors of \cite{Parlitz1990} for the first time demonstrated the emergence of chaotic oscillations of a single bubble through the cascade of period doubling bifurcations, which was later confirmed for an encapsulated bubbles in \cite{Macdonald2006,Carroll2013,Dzaharudin2013}. Furthermore, oscillations of coupled bubbles are more complicated than those of an isolated bubble due to the increasing dimensionality of the phase space. There are also new phenomena like synchronization, spatial arrangement and size distributions that appear when one considers even small clusters of bubbles \cite{Dzaharudin2013,Ohl2019,Garashchuk2019,Garashchuk2020,Garashchuk2021,Sojahrood2021}.
For example, in works \cite{Garashchuk2019,Garashchuk2021} two new scenarios of the onset of hyperchaotic oscillations were reported in a model of two interacting ultrasound contrast agents, where the former scenario is based on the formation of the discrete Shilnikov attractor and the latter on a specific mechanism of the destruction of synchronous oscillations.

It is also worth noting that the dynamics of microbubbles is inherently multistable in a physically and medically relevant range of the parameters \cite{Parlitz1990,Garashchuk2018,Garashchuk2019}, and, thus, requires careful consideration during the bifurcation analysis (see, e.g. \cite{Prasad2018,Prasad2018a,Prasad2023}). For example, in \cite{Garashchuk2020} it was shown that multistability has an influence on the stability of the synchronous oscillations of two interacting gas bubbles.

In this work we study nonlinear dynamics of three encapsulated microbubbles interacting via the Bjerknes forces. To the best of our knowledge, this problem has only been partially considered in works \cite{Macdonald2006,Dzaharudin2013} where the Morgan model for the bubbles shell was used, which has been suggested to be incorrect (see discussion in \cite{Doinikov2011,Doinikov2013}). In addition, the multistability of bubbles oscillations was not taken into account in bifurcation analysis in the above mentioned works. Here we use the de Jong model \cite{deJong1992,Marmottant2005,Doinikov2013} to describe the thin lipid shell of bubbles. We also suppose that the spatial configuration of bubbles is completely symmetric, i.e. they are located at the vertices of an equilateral triangle and the bubbles have the same equilibrium radii. A more detailed study of the oscillations of three interacting bubbles is necessary to uncover new dynamical effects that occur in small clusters of contrast agents and have not been observed previously. Investigation of the mechanisms of their appearance can provide a better understanding of the behavior of ensembles of interacting contrast agents that are used in applications.

We begin with identifying possible types of bubbles oscillations in the considered setting. We show that dynamics of three interacting bubbles can be periodic, quasiperiodic, chaotic, hyperchaotic with two and three positive Lyapunov exponents (LEs), chaotic with an additional zero (indistinguishable from zero in numerical experiments) LE and hyperchaotic with two positive LEs and an additional zero LE. Since the mechanisms of emergence of quasiperiodic, chaotic and hyperchaotic (with two positive LEs) oscillations are clear \cite{Garashchuk2018, Garashchuk2019, Garashchuk2021}, in this paper we focus on the onset of hyperchaotic attractors with three positive LEs and chaotic and hyperchaotic attractors with an additional zero LE.

On the one hand, the above-mentioned attractors are completely new in the context of the bubbles dynamics and hyperchaotic attractors with an additional zero Lyapunov exponent are novel in principle. On the other hand, there is no general theory describing how such attractors arise. Only a few specific results have been obtained in this direction (see \cite{Stankevich2020zero, Grines2022, Kazakov2023, Soldatkin2024}). For instance, in \cite{Soldatkin2024} it was demonstrated  that three-dimensional instabilities can appear within an attractor as a result of a sequence of period-doubling bifurcations or alternating period-doubling and Neimark-Sacker bifurcations. As far as chaotic attractors possessing an additional zero LE are concerned, it is known that they can arise in dynamical systems under quasiperiodic forcing \cite{Stankevich2021} or as a result of a cascade of torus-doubling bifurcations \cite{Stankevich2020zero}.

In the absence of the general framework explaining the onset of the hyperchaotic attractors characterized by three positive LEs and strange attractors with an additional zero LE, our main aim is to investigate the mechanisms of their appearance. We propose several scenarios of the onset of these type of attractors. In the case of hyperchaotic attractors with three positive LEs we connect their appearance to either the pairs of successive period-doubling bifurcations or Neimark-Sacker bifurcations of unstable periodic orbits (with one-dimensional unstable manifold) that lie inside these attractors. We believe that attractors with an additional zero LE appear due to the partial synchronization of two bubbles, which are affected by the quasiperiodic signal generated by the third bubble oscillating under the influence of the external periodic force. We believe that these scenarios can be observed in other dynamical systems and systems of coupled oscillators (see, e.g. \cite{Stankevich2020zero,Stankevich2023,Kapitaniak2022,Kapitaniak2023}).

The rest of this work is organized as follows. In the next Section we introduce the governing dynamical system for the oscillations of three interacting encapsulated spherical bubbles in a viscous liquid. Section \ref{sec:dyn} is devoted to the general overview of possible types of dynamics in the studied model, while in Sections 4 and 5 hyperchaotic oscillations with three positive LEs and strange attractors with an additional zero LE are studied, respectively.

\section{Main system of equations}
\label{sec:eq}

We consider three gas bubbles encapsulated into viscoelastic shells in a viscous compressible liquid with pairwise interactions described by the Bjerknes force \cite{Takahira1995,Mettin1997,Ida2002,Doinikov2004,Dzaharudin2013,Macdonald2006}. The corresponding mathematical model has the form
\begin{equation}
\label{eq:2bbls}
  \rho \left(R_{i} \ddot{R_{i}} +\frac{3}{2}\dot{R}_{i}^{2}\right)=P_{i}-\sum\limits_{j\neq i}\frac{d}{dt}\left(\frac{R_{j}^{2}\dot{R}_{j}}{d_{ij}}\right),
\end{equation}
where
\begin{multline}
	\label{eq1a}
	P_{i}= \left(P_{0}+\frac {2\sigma}{ R_{i,0}} \right)  \left( \frac { R_{i,0}}{R_{i}  } \right) ^{3\kappa} \left( 1-\frac {3\kappa }{c} R_{i,t}
	\right) -\frac {4\mu R_{i,t}  }{R_{1}}-\frac {2\sigma}{R_{i}}-P_{0}-  \\
	-   a \sin \left( \omega\,t \right) -4 \chi \left( \frac{1} {R_{i,0}}- \frac{1}{R_{i}}\right) -\frac {4\kappa_{s} R_{i,t} }{ R_{i} ^{2}}, \quad i,j=1,2,3.
\end{multline}
Here $t$ is time, $R_i(t)$, $i=1,2,3$ are radii of each bubble, $d_{ij}$ is the distance between centers of bubbles ($d_{ij}=d_{ji}$), $P_0 = P_{\mathrm{stat}} - P_v$, where $P_{\mathrm{stat}}$ denotes the static pressure and $P_v$ -- is the vapor pressure, $a$ is the amplitude of the external pressure field with the cyclic frequency $\omega$, $\sigma$ is the surface tension, $\rho$ is the liquid's density and $\eta_L$ is its viscosity, $c$ is the speed of sound, $\gamma$ is the polytropic exponent, $\chi$ is the elasticity of the shell, $\kappa_s$ is its viscosity, and $R_{i,0}$ is the equilibrium radius of the i-th bubble.

In the system \eqref{eq:2bbls} we use the de Jong model \cite{deJong1992,Marmottant2005,Doinikov2013} to describe bubbles' shells (see the last two terms in \eqref{eq1a}). In order to take into account liquid's compressibility we use the approach proposed in \cite{Marmottant2005}.

Equation \eqref{eq:2bbls} can be converted into an equivalent seven-dimensional autonomous dynamical system with phase variables $R_{i}$, $ \dot R_{i}$, $i=1,2,3$ and $\theta$, where the last equation is $\dot \theta=1$. The corresponding initial condition is $\theta(0)=0$. In what follows we will use this dynamical system for all numerical calculations.

Natural control parameters in \eqref{eq:2bbls} are the amplitude of the ultrasound field $a$ and its frequency $\omega$. The distances between bubbles $d_{ij}$, $i\neq j$, $i,j=1,2,3$ can also be considered as control parameters since if we assume that our three bubbles are a part of a cluster, its density would regulate the interbubble distances. In order to reduce the dimension of the parameter space throughout this work we fix the value $\omega=2.87 \cdot 10^7 \cdot s^{-1}$, which corresponds to biomedical applications of microbubbles. We also assume that the bubbles have the same equilibrium radii ($R_{0,1}=R_{0,2}=R_{0,3}=R_{0}$).

To perform numerical calculations we use the nondimensional variables $r_i, \tau$: $R_{i}=R_{0}r_{i}$, $t=\omega_{0}^{-1}\tau$, where $\omega_{0}^{2}=3\kappa P_{0}/(\rho R_{0}^2)+2(3\kappa-1)\sigma/R_{0}+4\chi/R_{0}$ is the natural frequency of bubbles oscillations in the linear approximation. The nondimensional speeds of the bubbles boundary can be expressed as follows: $u_{i}=dr_{i}/d\tau=\dot{R}_{i}/ (R_{0}\omega_{0})$, $i=1,2,3$. We construct the Poincar\'e map as a cross section by taking the values of the phase variables at each period of the external pressure field: $\{r_i(\tau_k), u_i(\tau_k)\}, i=1,2,3$, where $\tau_k = (2 \pi k \omega_{0})/\omega $ and $k \in \mathbb{N}$. The points in the Poincar\'e map belong to a 6-dimensional space. Thus, for the visualization purpose we use suitable projections, e.g. on the $(r_1, r_2, r_3)$ or $(r_1, u_1, r_2)$ subspaces.

We consider the fully symmetrical spatial configuration of the bubbles, assuming that they are located in the vertices of an equilateral triangle. Therefore, there is one effective parameter $d=d_{12} = d_{23} = d_{13}$, which describes the distance between bubbles. This symmetric configuration leads to the existence of a discrete symmetry corresponding to the exchange of any pair of indices $i \longleftrightarrow j$, $i =1,2,3, j = 1,2,3, i \neq j$.  Thus, there exist three 5-dimensional hyperplanes of partial synchronization, described by the linear constraints $S_{12}: R_{1} = R_{2}, \dot{R}_1 = \dot{R}_2$,  $S_{13}: R_{1} = R_{3}, \dot{R}_1 = \dot{R}_3$ and  $S_{23}: R_{2} = R_{3}, \dot{R}_2 = \dot{R}_3$, in the 7-dimensional phase space of system \eqref{eq:2bbls},  and a 3-dimensional hyperplane of complete synchronization $S: R_{1} = R_{2} = R_3, \dot{R}_1 = \dot{R}_2 = \dot{R}_3$. Notice that $S$ can be considered as the intersection of any pair of the manifolds of partial synchronization: $S = S_{12} \cap S_{13} = S_{12} \cap S_{23} = S_{13} \cap S_{23}$. The trajectories belonging to $S_{ij}$ correspond to in-phase synchronous oscillations of $i$th and $j$th bubbles, while the trajectories lying within $S$ describe complete synchronization of the oscillations of all the bubbles. Note that we assume that a point of a trajectory belongs to the manifold of complete or partial synchronization if its distance in the $L_{1}$ norm from the corresponding hyperplane is less than $10^{-6}$

Considering completely synchronous oscillations, the governing system of equations can be reduced to a 3-dimensional one by substituting the corresponding linear constraints into \eqref{eq:2bbls}:
\begin{equation}
\label{eq:synchr}
	\rho \left(R_{1} \ddot{R_{1}} +\frac{3}{2}\dot{R}_{1}^{2}\right)=P_{1}-2\frac{d}{dt}\left(\frac{R_{1}^{2}\dot{R}_{1}}{d}\right).
\end{equation}
Notice that the difference from the model of an isolated bubble comes from the last term $2\frac{d}{dt}\left(\frac{R_{1}^{2}\dot{R}_{1}}{d}\right)$. Moreover,  completely synchronous oscillations of two interacting contrast agents can also be described by \eqref{eq:synchr} but without the multiplier $2$ in the last term. The last term from \eqref{eq:synchr} can be rewritten as $\frac{d}{dt}\left(\frac{R_{1}^{2}\dot{R}_{1}}{(d/2)}\right)$. This means that completely synchronous oscillations, which appear in the model of two interacting bubbles, can also be observed in the case of three interacting contrast agents equidistant from each other. This requires the same values of the parameters $a$ and $\omega$ and two times bigger distance between bubbles.

\section{Nonlinear dynamics of three interacting microbubbles}
\label{sec:dyn}

In this section we determine various types of dynamics of three interacting microbubble contrast agents governed by system \eqref{eq1a}. We consider a completely symmetrical spatial configuration of the bubbles located at the vertices of an equilateral triangle. We use the magnitude of the external pressure field and the distance between bubbles as the control parameters. We determine possible types of bubbles oscillations in a wide range of the control parameters.

\begin{figure}[ht]
	\centering{
		 \includegraphics[width=0.9\linewidth]{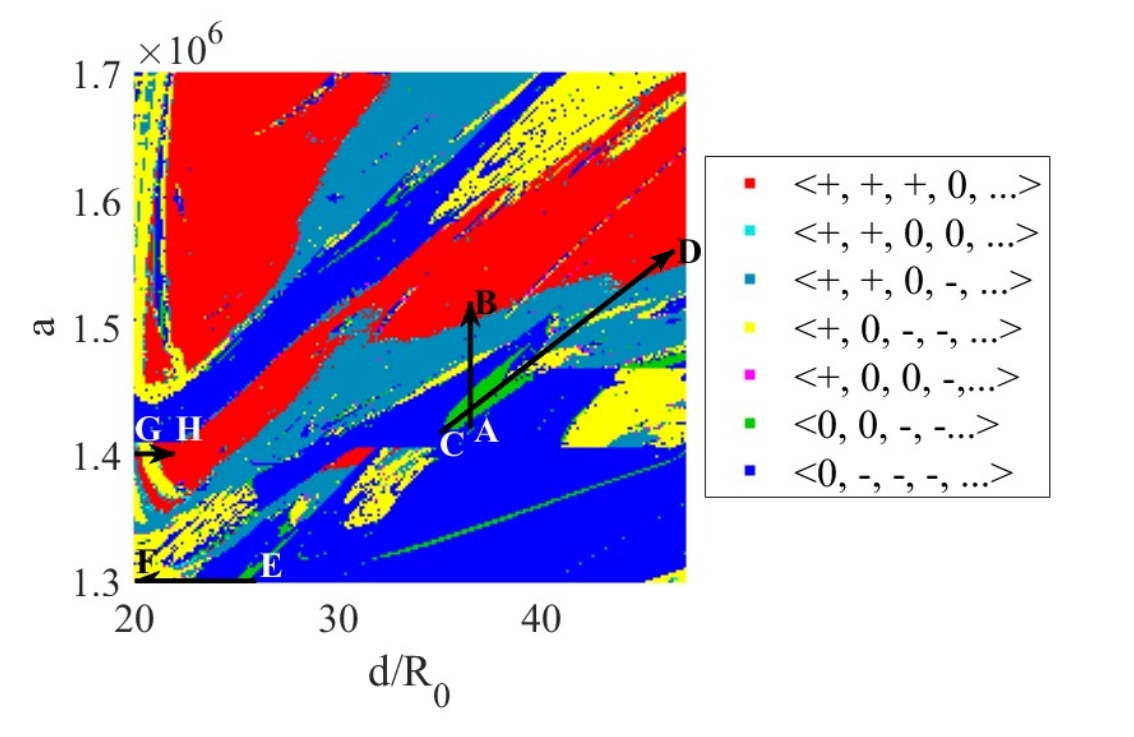}

	}
	\caption{The chart of LEs on the parameter plane $(d/R_{0}, a)$: $1.3 < a/10^6 < 1.7$, $20 < d/R_{0} < 47$. The detailed one-parametric analysis is carried out along the highlighted routes AB, CD, EF, and GH.
	}

	\label{fig:chart}
\end{figure}

In Fig. \ref{fig:chart} we demonstrate the chart of types of oscillations on the parameter plane $(d/R_0, a)$. At each pixel of this chart we compute LEs using the standard procedure suggested in \cite{Benettin1980}\footnote{Before computing LEs we skip a sufficiently long transient process to be sure that an orbit reaches an attractor. Notice also that the dynamics in the system under consideration is often multistable. To suppress undesirable leaps between coexisting attractors we take initial conditions corresponding to a quasiperiodic regime at $(d/R_0, a/10^6)=(47,1.48)$, continue it along $a$-axis in both directions on the line $d/R_0=47, 1.3 < a/10^6 < 1.7$ and then we continue the resulting regimes from this line along the $d/R_{0}$-axis in the direction of decreasing $d/R_0$.}. Let us remark that any attractor of seven-dimensional system \eqref{eq:2bbls} always has at least one zero LE $\lambda_z=0$, corresponding to the translations along the attractor. The sings of the remaining six exponents $\lambda_1 \geq \lambda_2 \geq ... \geq \lambda_6$ determine the type of dynamics of the regimes found, see the color palette in Table~\ref{table:1}. We estimate that error in the numerical computations of LEs is less than $10^{-4}$ and consider an exponent to be equal to zero if its absolute value is less than $10^{-4}$.

\begin{table}[h]
	\begin{center}
	\begin{tabular}{|m{1cm}| m{8cm}| m{3cm}|}
		\hline
		Color & Dynamical regime & Signature of the spectrum of the LEs \\ [0.5ex] 
		\hline
		\begin{tikzpicture}
			\draw[fill=red,draw=black, thin] (0,0) rectangle (0.5,0.25);
		\end{tikzpicture}  & Hyperchaotic with three positive LEs & $<+, +, +, 0, -, ...>$ \\
		\begin{tikzpicture}
			\draw[fill=cyan,draw=black, thin] (0,0) rectangle (0.5,0.25);
		\end{tikzpicture} & Hyperchaotic with two positive LEs and one additional zero LE & $<+, +, 0, 0, -, ...>$ \\
		\begin{tikzpicture}
			\draw[fill=teal, draw=black,thin] (0,0) rectangle (0.5,0.25);
		\end{tikzpicture} & Hyperchaotic dynamics with two positive Lyapunov exponents & $<+, +, 0, -, ...>$ \\
		\begin{tikzpicture}
			\draw[fill=yellow,draw=black, thin] (0,0) rectangle (0.5,0.25);
		\end{tikzpicture} & Chaotic & $<+, 0, -, ... >$ \\
		\begin{tikzpicture}
			\draw[fill=magenta,draw=black, thin] (0,0) rectangle (0.5,0.25);
		\end{tikzpicture} & Chaotic with an additional zero LE & $<+, 0, 0, -, ... >$  \\
		
		\begin{tikzpicture}
			\draw[fill=green!80!black,draw=black, thin] (0,0) rectangle (0.5,0.25);
		\end{tikzpicture} & Quasiperiodic dynamics& $<0, 0, -, -, ... >$ \\
		
		\begin{tikzpicture}
			\draw[fill=blue!,draw=black, thin] (0,0) rectangle (0.5,0.25);
		\end{tikzpicture} & Periodic dynamics& $<0, -, -, ... >$ \\ [0.5ex]
		\hline
	\end{tabular}
	\end{center}
	\caption{Correspondence between the colors used in Fig.~\ref{fig:chart} and types of dynamics, which are defined according to the spectrum of LEs.}
	\label{table:1}
\end{table}

From Fig.~\ref{fig:chart} one can see that there are at least five types of strange attractors in system \eqref{eq:2bbls}: simple chaotic attractors (with one positive and one zero LEs), chaotic attractors with an additional zero LE, simple hyperchaotic attractors (with two positive and one zero LEs), hyperchaotic attractors with three positive LEs and hyperchaotic attractors with two positive and two zero LEs.

Simple chaotic and hyperchaotic attractors for systems describing the dynamics of gas bubbles and the bifurcation mechanisms of their appearance were studied in details \cite{Garashchuk2018, Garashchuk2019}. In particular it was shown that synchronous chaotic dynamics arise via the cascade of period-doubling bifurcations \cite{Feigenbaum}, while asynchronous chaotic attractors appear via the Afraimovich--Shilnikov scenario of the destruction of an invariant torus \cite{AfrShil1983} (see also \cite{broer1998towards}). Furthermore, the onset of asynchronous hyperchaotic dynamics is governed by the Shilnikov scenario \cite{GGS12} leading to the attractor containing a saddle-focus periodic orbit with two-dimensional unstable manifold. Since these scenarios have been extensively studied earlier \cite{Garashchuk2018, Garashchuk2019, Garashchuk2020, Kazakov2023}, we do not discuss them in this work. Instead, we focus on the hyperchaotic attractors with three positive LEs, hyperchaotic attractors with two positive and an additional zero LEs and chaotic attractors with an additional zero LE and the mechanisms of their appearance, because such regimes have not been previously observed in the models of interacting microbubbles.

\section{Onset of hyperchaotic dynamics with three positive Lyapunov exponents}
\label{sec:hyper}

In this section we study bifurcation scenarios leading to the onset of hyperchaotic dynamics with three positive LEs. We consider two paths AB and CD (see Fig. \ref{fig:chart}) leading to hyperchaotic oscillations. Let us remark that below, we use indices $(n,m)$ to denote a periodic orbit on the Poincar\'e map, which means that this orbit has $n$-dimensional stable and $m$-dimensional unstable invariant manifolds. Notice that in all cases $m+n=6$.

\begin{figure}[ht]
	\centering{
		
		 \includegraphics[width=\linewidth]{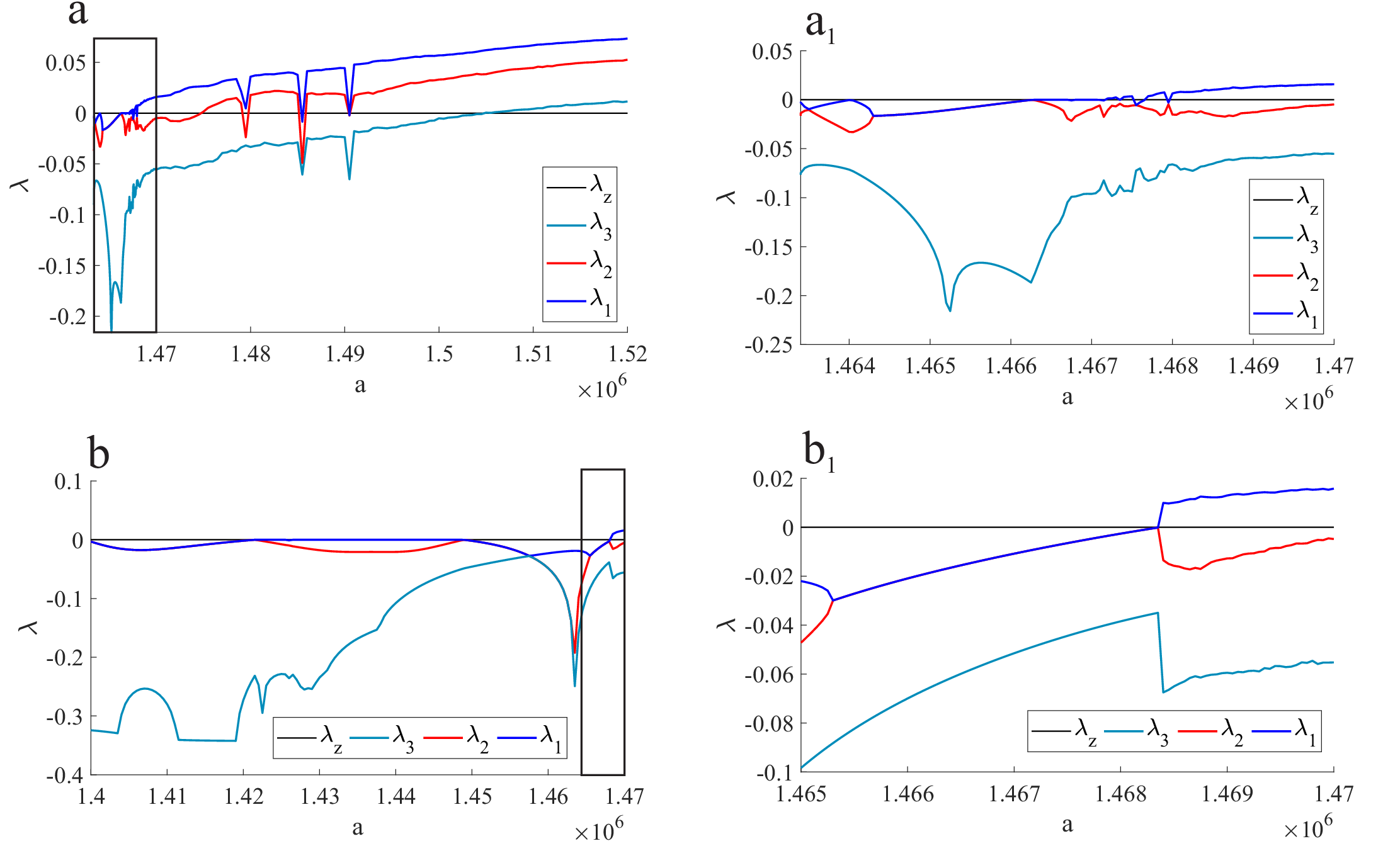}
	}
	\caption{Plots of the largest LEs along the path AB: (a) continuation of the period-8 orbit transforming into the Lorenz-like quasiattractor and then into the hyperchaotic regime with three positive LEs in $1.4634 < a/10^6 < 1.52$, enlarged image showing the emergence of chaos in $1.4634 < a/10^6 < 1.47$; (b) continuation of another regime existing on this route from $a/10^6=1.4$ to $a/10^6 = 1.47$, the enlarged fragment demonstrates the subcritical Neimark-Sacker bifurcation, accompanied by the transition to the chaotic attractor, $1.465 < a/10^6 < 1.47$.}
	\label{fig:AB}
\end{figure}

\begin{figure}[ht]
	\centering{
		\includegraphics[width=1.0\linewidth]{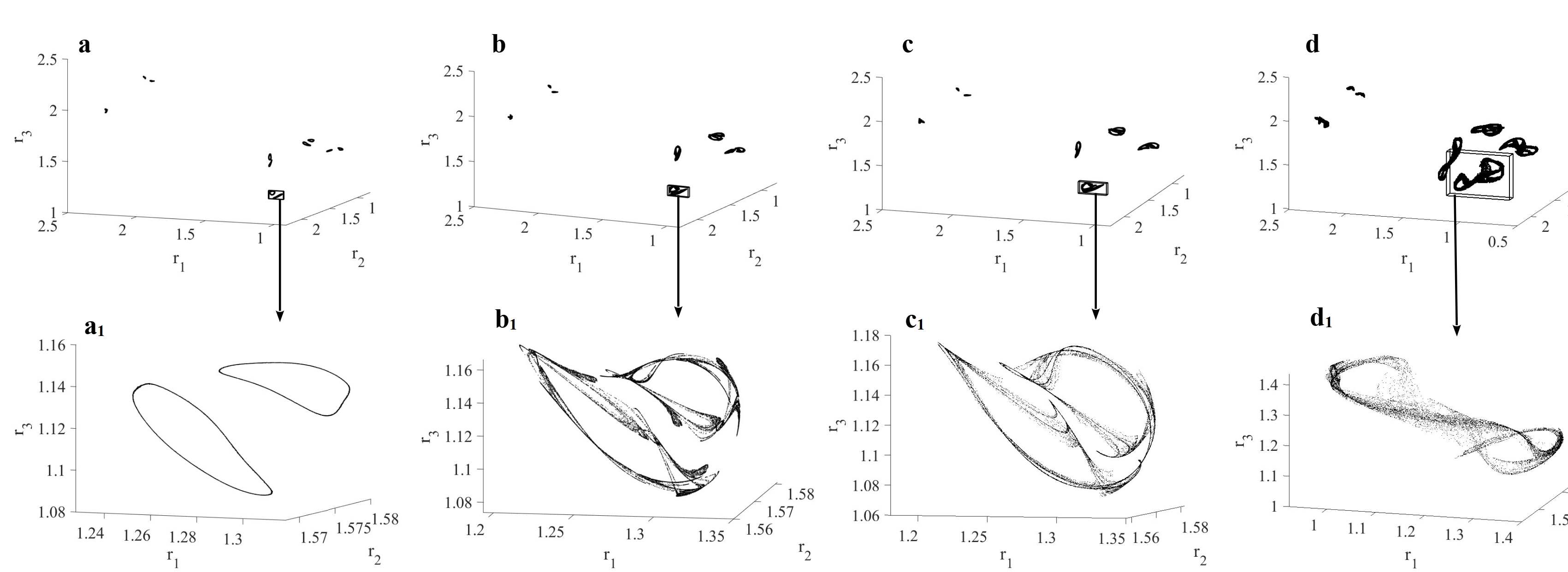}
	}
	\caption{Projections of the Poincar\'e maps of several attractors along the path AB on the subspace $(r_1, r_2, r_3)$: (a) the invariant torus with sixteen components at $(d/R_0, a/10^6) = (36.5 ,1.4667)$; (b) the chaotic attractor with sixteen components at $(d/R_0, a/10^6) = (36.5 ,1.4675)$; (c) the Lorenz-like chaotic attractor at $(d/R_0, a/10^6) = (36.5 ,1.4679)$; (d) the chaotic attractor at $(d/R_0, a/10^6) = (36.5 ,1.47)$.}
	\label{fig:AB_sec}
\end{figure}

First, we consider the path AB, which corresponds to $d/R_0 = 36.5$ and $a/10^{6}\in[1.4,1.52] $. The plots of the four largest LEs for the continuations of two attractors existing along this path are shown in Fig. \ref{fig:AB}. Several typical phase portraits on the Poincar\'e map are shown in Fig.~\ref{fig:AB_sec}. At the point $(d/R_0, a/10^6) = (36.5, 1.4634)$ a saddle-node bifurcation occurs, giving rise to an orbit of period eight with respect to the driving force.

If to increases $a$, one can see a period-doubling bifurcation, and, thus, a stable period-16 orbit appears. Later this orbit undergoes the supercritical Neimark--Sacker bifurcation, and we observe a quasiperiodic regime (stable invariant torus) with sixteen components (see Fig. \ref{fig:AB_sec}a). Further along the path, this torus loses its smoothness. As a result, the torus-chaos attractor with sixteen components and one positive LE arises (see Fig. \ref{fig:AB_sec}b).

A further increase in $a$ leads to a pairwise merging of sixteen components of the attractor into eight. This merging is due to the collision of attractor components with the stable manifold of the period-8 saddle orbit. This orbit has the five-dimensional stable and one-dimensional unstable invariant manifolds and lies between a pair of components of the sixteen-component attractor. Before the merger, left (right) piece of its unstable manifold goes to the left (right) component of the attractor, while the stable manifold separates these two components (Fig.~\ref{fig:AB_sec}b$_1$). Thus, eight triplets are formed from the components of the saddle orbit and the corresponding pairs of the chaotic attractor components. Each of these triplets merge into a new larger component during the collision, thus, forming an eight-component chaotic attractor (see Fig. \ref{fig:AB_sec}c).

The regime observed after the collision is called the discrete Lorenz-like attractor \cite{Gonchenko2005}. This is a generalization of the classical Lorenz-like attractors \cite{Lor63, guckenheimer1979structural, williams1979structure, ABS77, ABS82, shil1993bifurcations, shil1993normal} to the case of systems with discrete time. The main feature of such attractors is that they can be pseudohyperbolic \cite{Gonchenko2021, Kaz2024}. However, in our case the saddle value of the period-8 saddle orbit is less than one in absolute value, therefore here we observe a Lorenz-like quasiattractor \cite{evgenii2021lorenz}. If we increase $a$ further, this attractor expands and collides with another 8-component Lorenz-like attractor coexisting with it before the collision. Due to the symmetry of model \eqref{eq:2bbls}, the 8-component Lorenz-like attractor has a symmetric analogue. Their absorbing domains are separated by the stable invariant manifold of a symmetric saddle orbit which lies between this pair of attractors. At some point these attractors touch this stable manifold. After this boundary crisis we observe the attractor shown in Fig.~\ref{fig:AB_sec}d.

\begin{figure}[ht]
	\centering{
		 \includegraphics[width=0.8\linewidth]{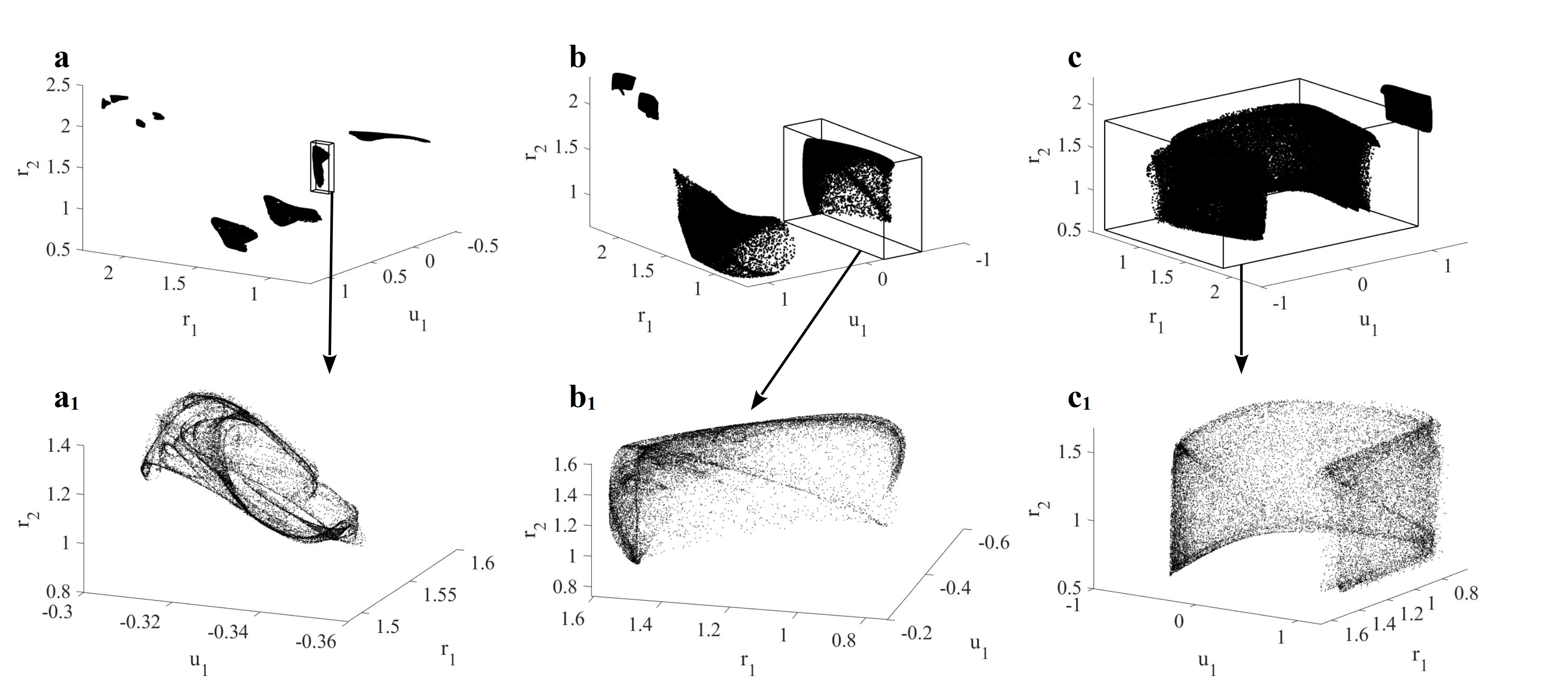}
	}
	\caption{Projections of the Poincar\'e maps of three attractors along the path AB on the subspace $(r_1, u_1, r_2)$: (a) the chaotic attractor with one positive LE at $(d/R_0, a/10^6) = (36.5 ,1.472)$, (b) the hyperchaotic attractor with two positive LEs at $(d/R_0, a/10^6) = (36.5 ,1.484)$ and (c) the hyperchaotic attractor with three positive LEs at $(d/R_0, a/10^6) = (36.5 ,1.52)$.
	}
	\label{fig:AB_sec_2}
\end{figure}

Increasing $a$ along the path $AB$, we see that the chaotic attractor expands in the phase space. Consequently, a smooth onset of hyperchaotic dynamics with two, and then three, positive LEs is observed at $a/10^6 \approx 1.475 $ and $a/10^6 \approx 1.505$, respectively (see Fig. \ref{fig:AB}a). Different components of the attractor begin to overlap visually in the projections onto the subspace $(r_1, r_2, r_3)$. Therefore, we use the coordinates $(r_1, u_1, r_2)$ to show further transformations of the attractor. In Fig.~\ref{fig:AB_sec_2}a we show an eight-component chaotic attractor with one positive LE. A hyperchaotic attractor characterized by two positive LEs is shown in Fig. \ref{fig:AB_sec_2}b. This attractor consists of four components, it appears after the pairwise merger of eight components into four. In Fig.~\ref{fig:AB_sec_2}c we show a hyperchaotic attractor with three positive LE consisting of only two components.

Based on the ideas presented in \cite{Kazakov2023} and \cite{Soldatkin2024} (see also \cite{Kapitaniak2022}), we believe that this attractor becomes a hyperchaotic one with three positive LEs according to one of the two following mechanisms. In the first one, periodic saddles of type (5,1) that separate different components of the attractor undergo two successive period-doubling cascades transforming them, first, into saddles of type (4,2) and, then, into saddles (3,3). According to the second one, periodic saddles of type (5,1) undergo the cascade of supercritical Neimark--Sacker \cite{sataev2021cascade} bifurcations which gives rise to saddles of type (3,3). After the merger of components of the attractor, these saddle periodic orbits with three-dimensional unstable manifold become a part of the attractor. As a result, we observe three positive LEs in numerical experiments.

Due to the cumbersome nature of the model under consideration we cannot trace the saddle periodic orbits and determine which particular scenario takes place along the one-parametric paths discussed here. However, the visual similarity of the transformations of the phase portraits and the smooth growth of the LEs to the numerical characteristics similar to those obtained in \cite{Soldatkin2024} gives us some confirmation of the described scenarios.

Notice also that there is multistability along the route AB (see Fig. \ref{fig:AB}b). There exists a stable orbit of period eight at $(d/R_0, a/10^6) = (36.5, 1.4)$. Later, it undergoes the supercritical Neimark--Sacker bifurcation, leading to a quasiperiodic regime with eight components. After that an inverse Neimark--Sacker bifurcation occurs, making the 8-periodic orbit stable again. Further increasing $a$, we observe a subcritical Neimark--Sacker bifurcation, making this periodic orbit unstable (see Fig. \ref{fig:AB}b). After it happens, we observe only the Lorenz-like chaotic attractor, discussed above.

 \begin{figure}[!ht]
	\centering{
		 \includegraphics[width=\linewidth]{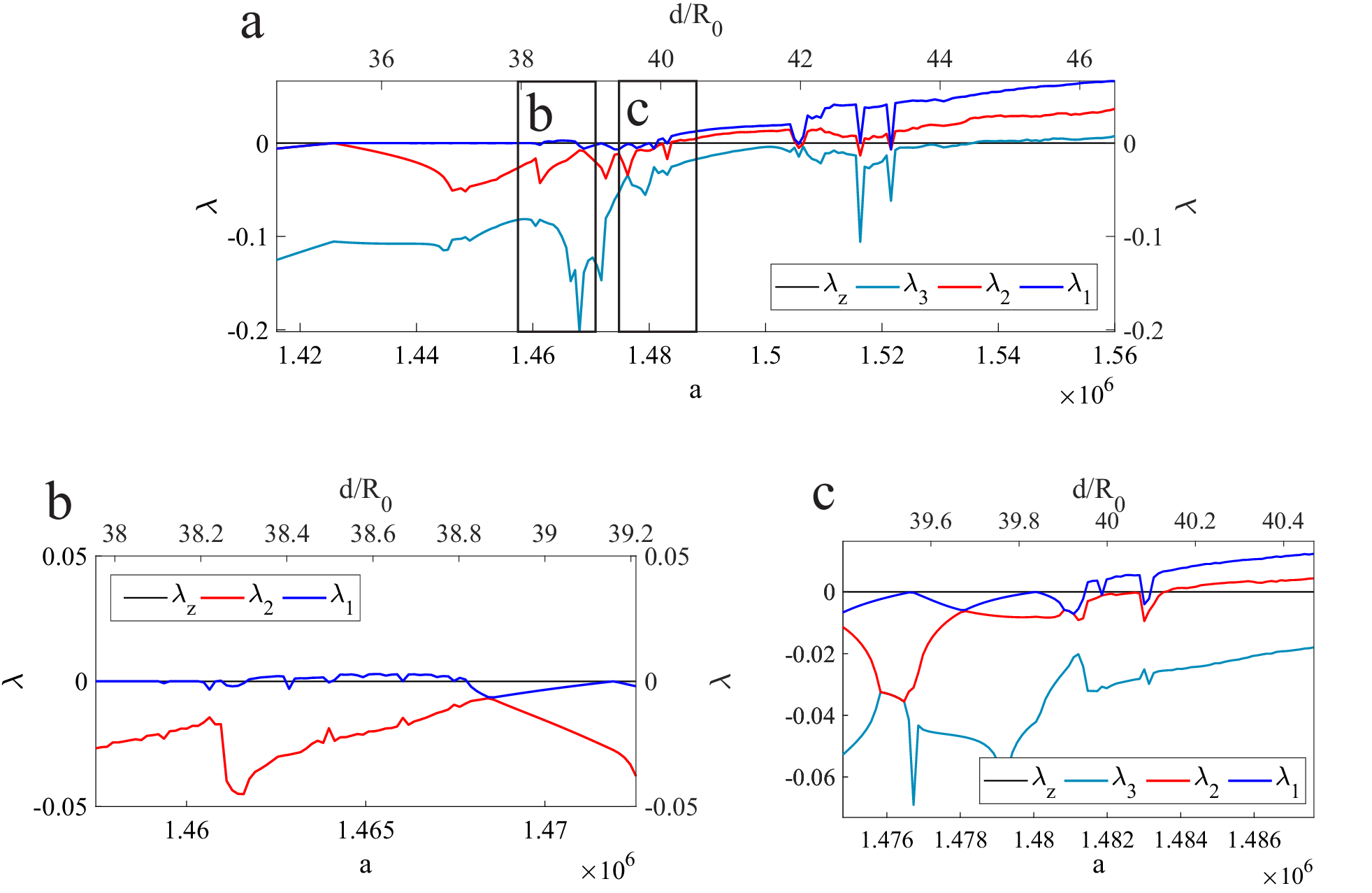}
	}
	\caption{The dependence of the largest LEs along the path CD on the parameter $a$: (a) plots of the four largest LEs for the path CD; (b) enlarged fragment showing the onset of the torus-chaos attractor; (c) enlarged fragment showing the transition to a hyperchaotic attractor with three positive LEs.}
	\label{fig:CD}
\end{figure}

\begin{figure}[!hb]
	\centering{
		 \includegraphics[width=\linewidth]{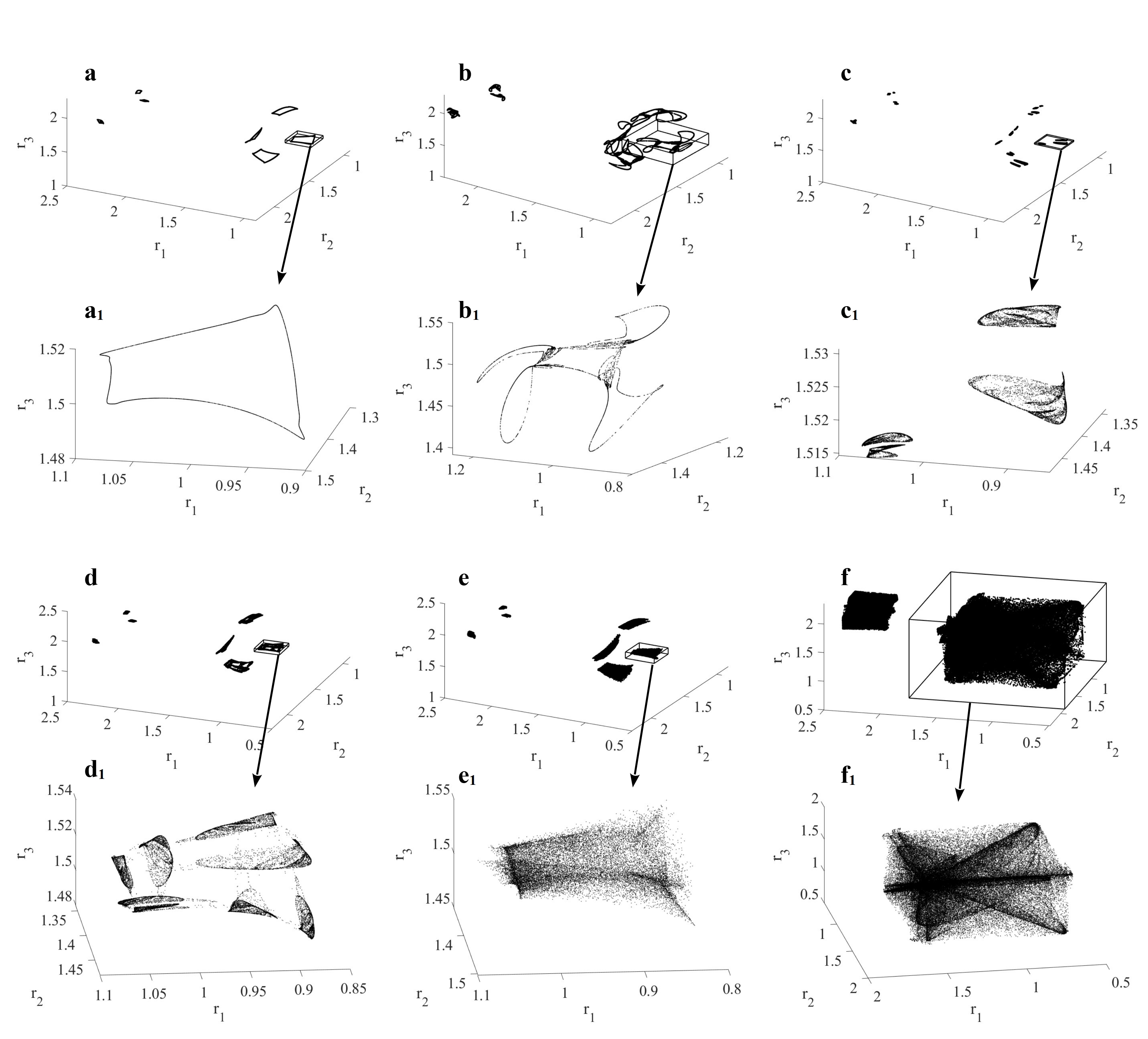}
	}
	\caption{Projections of the Poincar\'e maps for several attractors along the path CD on the subspace $(r_1, r_2, r_3)$ . (a) $(d/R_0, a/10^6) = (38, 1.4579)$, (b) $(d/R_0, a/10^6) = (38.6, 1.4652)$, (c) $(d/R_0, a/10^6) = (39.9666, 1.4816)$, (d) $(d/R_0, a/10^6) = (40.03, 1.4823)$, (e) $(d/R_0, a/10^6) = (41.6, 1.5012)$, (f) $(d/R_0, a/10^6) = (46.5, 1.56)$.	}
	\label{fig:CD_sec}
\end{figure}

Another route where we observe the onset of hyperchaotic attractors with three positive LEs is the path CD. This one-parametric route is taking place on the line between two points, $C: (d/R_0, a) = (35, 1.1416), D: (d/R_0, a) = (46.5, 1.56)$. The plot of the maximal LEs along this path is shown in Fig.~\ref{fig:CD}.  At the point C there is a stable limit cycle. The first bifurcation along this path is again the supercritical Neimark-Sacker bifurcation, that results in a stable invariant torus (see Fig.~\ref{fig:CD_sec}a). Moving further along this path we observe the Afraimovich--Shilnikov scenario \cite{AfrShil1983} of the torus destruction. As a result, we obtain the so-called torus-chaos attractor (see the graph of the maximal LEs in Fig.~\ref{fig:CD}b and the corresponding phase portrait on the Poincar\'e map in Fig.~\ref{fig:CD_sec}b). Then, one can see a window of stability, where a stable periodic orbit appears. This orbit undergoes three period-doubling bifurcations within this window of stability. Some of these bifurcations can be seen in the right part of Fig. \ref{fig:CD}b when $\lambda_1$ approaches zero and the rest of them are shown in Fig. \ref{fig:CD}c. Then we observe a sudden transition to a chaotic attractor shown in Fig.~\ref{fig:CD_sec}c. Here one can see four disconnected components inside each part of the phase space, where one part of the eight-component attractor used to exist. These components begin to merge while the chaotic attractor continues to expand along the path CD (see Fig.~\ref{fig:CD_sec}d). Notice that hyperchaotic dynamics appear after the merger of the smaller components (see the hyperchaotic attractor with two positive LEs in Fig.~\ref{fig:CD_sec}e). Then, this process continues with the two groups consisting of four components of the attractor. The moment when only two large components remain is associated with the onset of hyperchaotic dynamics with three positive LEs (see Fig.~\ref{fig:CD_sec}f). Apparently, the underlying bifurcation mechanisms of the transition from chaotic dynamics with one positive LE to hyperchaotic dynamics with three positive Lyapynov exponents are the same as discussed earlier for the path AB.

\section{Attractors with an additional zero Lyapunov exponent}
\label{sec:addzero}

\begin{figure}[ht]
	\centering{
		 \includegraphics[width=\linewidth]{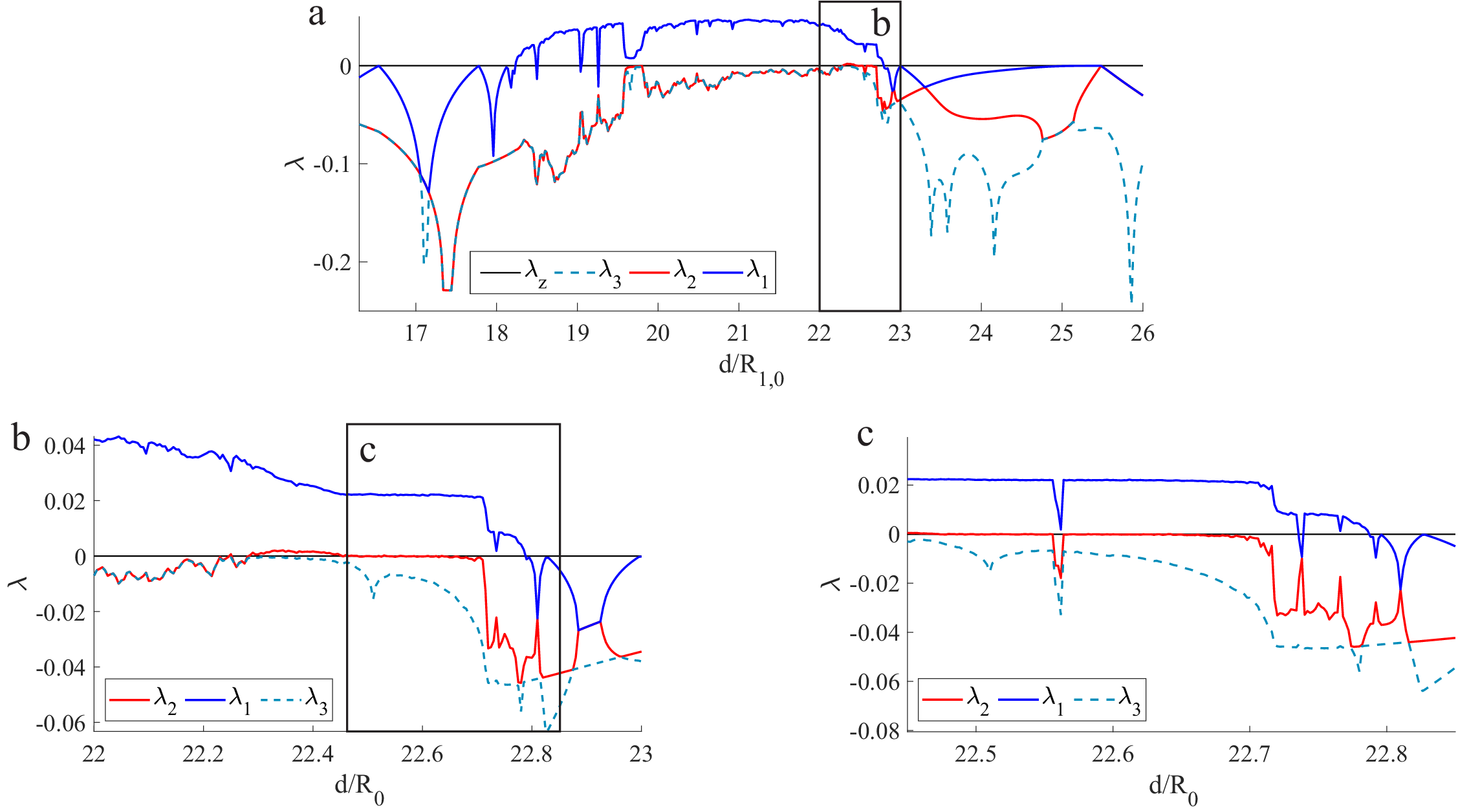}
	}
	\caption{(a) The dependence of the four largest LEs on the parameter $d/R_{0}$ for the path EF; (b), (c) enlarged fragments of this plot for $22 < d/R_{0} < 23$ and $22.45 < d/R_{0} < 22.85$ respectively, demonstrating the existence of a chaotic attractor with an additional zero LE. }
	\label{fig:EF}
\end{figure}

In this section we discuss the existence of chaotic and hyperchaotic attractors with an additional zero LE in system \eqref{eq:2bbls}. Up to now, chaotic attractors with an additional zero LE have been reported in a number of systems (both with continuous and discrete time) from various applications (see, e.g.~\cite{Broer2002, Broer2005, Gonchenko2005, Broer2010, Stankevich2020zero, Karatetskaia2021, Kazakov2023}). Nevertheless, in many cases it is not known how such attractors arise and why they exist in sufficiently large domains of the parameter space.

Explanations of this phenomenon are given only in some special cases: first, when the corresponding dynamical system is under the influence of an external quasiperiodic perturbation~\cite{Broer2010}; second, when the system, in fact, has a smaller effective dimensionality due to the existence of additional first integrals \cite{BorKazSat2016}; third, when such strange attractors are found near some codimension-3 bifurcations \cite{Gonchenko2005, Karatetskaia2021, Grines2022, Kazakov2023}. To the best of our knowledge, systems with hyperchaotic attractors characterized by additional zero LEs have not been reported previously.

\begin{figure}[ht]
	\centering{
		 \includegraphics[width=\linewidth]{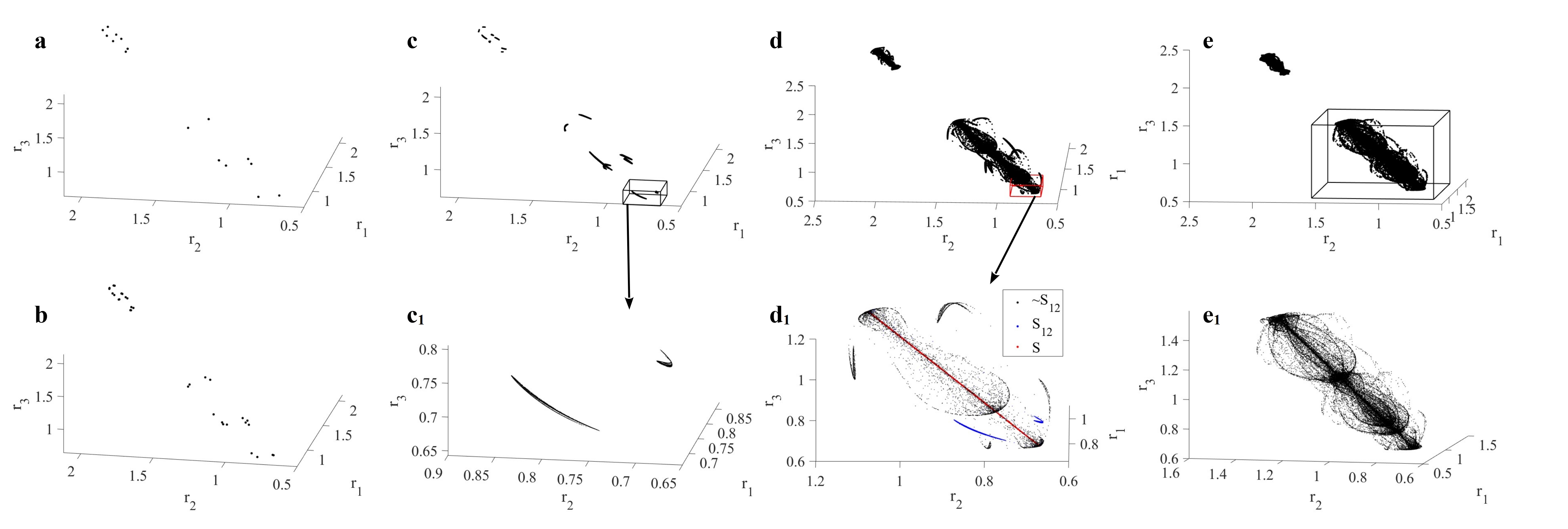}
	}
	\caption{Projections of the Poincar\'e maps of several attractors along the path EF on the subspace $(r_1, r_2, r_3)$: (a) partially synchronous limit cycle at $(d/R_0, a/10^6) = (23.0 ,1.3)$; (b) partially synchronous periodic regime (after several period-doubling bifurcations) at  $(d/R_0, a/10^6) = (22.84 ,1.3)$; (c) partially synchronous attractor at $(d/R_0, a/10^6) = (22.722 ,1.3)$ and a magnified image of one of its components; (d) chaotic attractor at  $(d/R_0, a/10^6) = (22.7712 ,1.3)$, where blue dots in the enlarged part correspond to the points in $S_{12}$ (partial synchronization), red dots -- to the points in $S$ (complete synchronization), and black to all the rest, here we use a red cube instead of the black one simply to make it visible on the background of the black dots; (e) chaotic attractor with an additional zero LE at $(d/R_0, a/10^6) = (22.6 ,1.3)$.}
	\label{fig:EF_sec}
\end{figure}

We begin with the path EF: $a/10^6 = 1.3, d/R_0 \in [16.5, 26]$, going down from high to low values of $d/R_0$ and trace the bifurcations that occurre with the stable period-8 orbit and the attractors that are born from it. The plot of the four largest LEs is presented in Fig.~\ref{fig:EF}. Note that the period-8 orbit corresponds to partial synchronization of the bubbles oscillations: the phase trajectory belongs to $S_{12}$ ($r_1=r_2$, $u_1=u_2$), which means that the first and the second bubbles oscillate in-phase. As $d/R_0$ decreases, this periodic orbit undergoes a cascade of period-doubling bifurcations (see the phase portraits before and after several period-doubling bifurcations in Figs.~\ref{fig:EF_sec}a and ~\ref{fig:EF_sec}b), resulting in a partially synchronous chaotic attractor shown in Fig.~\ref{fig:EF_sec}c. One can see that this chaotic regime has sixteen components.

If we further decrease $d/R_0$, two processes happen simultaneously: the components of the chaotic attractor begin to merge and the partial synchronization is destroyed (see Fig.~\ref{fig:EF_sec}d). Notice that the inclusion of the asynchronous orbits into the attractor happens suddenly, simultaneously with the inclusion of completely synchronous orbits (see Fig. \ref{fig:EF_sec}d$_1$, where we use blue color for the points belonging to $S_{12}$, red color for the points in $S$, and black color for all the other points\footnote{Some of the black points belong to the sets $S_{13}$ or $S_{23}$, but we do not use separate colors for them in order not to overcomplicate the presentation}).

It is important to note that the blue set in Fig. \ref{fig:EF_sec}d$_1$ looks similar to the corresponding components of the attractor shown in Fig.~\ref{fig:EF_sec}c$_1$. Therefore, the completely synchronous chaotic saddle set and the asynchronous chaotic saddle set develop separately, while they are not attractive. The inclusion occurs due to the collision of the partially synchronous chaotic attractor with these sets. In our case, together with the destruction of the synchronization, the attractor gets an additional zero LE, which exists in a substantial interval of the parameter $d/R_0$ (see the plot of the maximal LEs in Fig.~\ref{fig:EF}c). The phase portrait of this attractor is shown in Fig.~\ref{fig:EF_sec}e.

System \eqref{eq:2bbls} is not subjected to the action of a quasiperiodic forcing, it has no first integrals. Moreover, we do not find any specific bifurcations of co-dimension two or three. Therefore, the reason for the appearance of a chaotic attractor with an additional zero LE is be unclear. However, one can hypothesise that the periodic external driving force together with the impact of the third bubble can create an influence on the cluster of two synchronously oscillating bubbles similar to an external quasiperiodic perturbation, and, hence, leading to appearance of an additional neutral direction.

\begin{figure}[ht]
	\centering{
		\includegraphics[width=\linewidth]{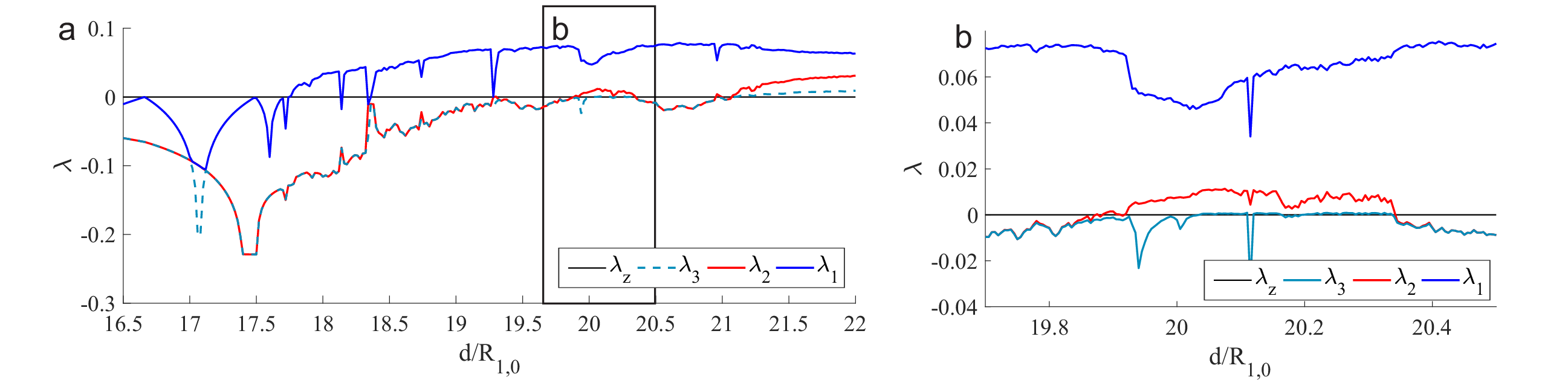}
	}
	\caption{(a) Graph of four largest LEs for the path GH; (b) its enlarged fragment for the interval $19.7 < d/R_{0} < 20.5$.
	}
	\label{fig:GH}
\end{figure}

\begin{figure}[ht]
	\centering{
		\includegraphics[width=\linewidth]{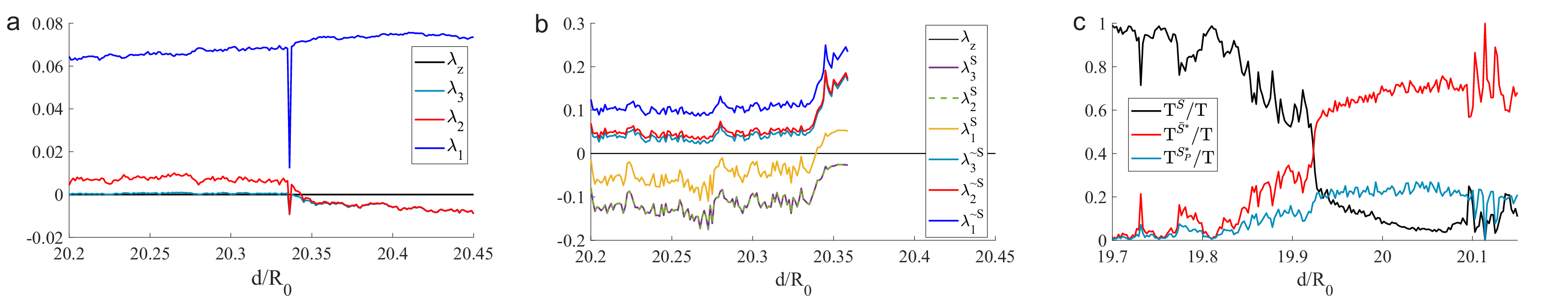}
	}
	\caption{(a) Graph of the largest LEs for the fragment of the path GH: $19.7 < d/R_{0} < 20.15$, (b) the largest LEs, corresponding to the synchronous and asynchronous phases of the the trajectory (when the trajectory spends sufficient time in each phase), (c) the ratios of time the associated with the synchronous, the asynchronous and the partially synchronous components of the attractor.
	}
	\label{fig:GH_split}
\end{figure}

Finally, we study the dynamics along the path GH: $a/10^6 = 1.4, d/R_0 \in [16.5, 22]$. The corresponding graph of the four largest LEs is presented in Fig.~\ref{fig:GH}. At the point G we observe a stable completely synchronous periodic orbit. Note that complete synchronization leads to the coincidence of the second and the third LEs. This is due to the symmetry of the effectively 3-dimensional attractor embedded into $S$. With increasing $d/R_{0}$ this periodic orbit undergoes a cascade of period-doubling bifurcations giving rise to a completely synchronous chaotic attractor.

Further down this path, we observe the destruction of synchronization via the bubbling transition mechanism. The same scenario was previously studied in the model of two interacting gas bubbles \cite{Garashchuk2021}. Here we use the notation $S^*_{ij} = S_{ij} \setminus  S$ to describe the points corresponding to partial synchronization (between bubbles $i$ and $j$), but not to complete synchronization of oscillations of all three bubbles\footnote{Any point corresponding to complete synchronization also belongs to any of the manifolds of partial synchronization: $\forall i,j: i \neq j \Rightarrow S \subset S_{ij}$.} Set $\bar{S}$ contains all the points outside of the hyperplane of complete synchronization. However, it is useful to consider also the points in $\bar{S}^* = \overline{S_{12}} \cap \overline{S_{13}} \cap \overline{S_{23}}$ that lie not only outside of $S$, but also do not belong to any of the manifolds of partial synchronization. For a typical trajectory on the attractor that we were following for the time $T$ we will denote the amount of time it spent inside the respective manifolds as $T^S$, $T^{S_{ij}}$, $T^{S^{*}_{ij}}$, $T^{\bar{S}}$, $T^{\bar{S}^{*}}$. We will also use the values $T^{S_P^*}= T^{S^*_{12}} + T^{S^*_{13}} + T^{S^*_{23}}$, which describe the total amount of time that the trajectory spends inside the manifolds of partial synchronization, excluding the states of complete synchronization.

Following the approach proposed in \cite{Garashchuk2021}, we present the graphs of the largest LEs for the interval $19.7 < d/R_0 < 20.15$ in Fig.~\ref{fig:GH_split}a and the graphs of the largest LEs, associated with the parts of the attractor lying in $S$ and in $\bar{S}$ separately in Fig.~\ref{fig:GH_split}b. Notice how the second and the third LEs, corresponding to the synchronous part of the attractor, are still equal after the appearance of the asynchronous component (see Fig.~\ref{fig:GH_split}b), due to the same symmetry considerations as for the completely synchronous attractor. We also show the fractions of time, associated with synchronous ($T^S/T$), partially synchronous ($T^{S_P^*}/T$) and completely asynchronous ($T^{\bar{S}^*}/T$) oscillations for a typical trajectory on the attractor in Fig.~\ref{fig:GH_split}c. The sum of these values is equal to one:  $T^S/T + T^{S_P^*}/T + T^{\bar{S}^*}/T = 1$. The transition from synchronous chaos to hyperchaos and then to hyperchaos with and additional zero LE happens gradually, as the partially synchronous and completely asynchronous components of the attractor increase, while the completely synchronous one decreases (see Fig.~\ref{fig:GH_split}). The appearance of the partially synchronous part of the attractor during this process suggests that some of the stable/unstable manifolds of the saddle orbits inside $S$ are directed towards $S_{ij}$, while $S_{ij}$ themselves have significant transversally stable domains. After the appearance of the hyperchaotic regime, consisting of the completely synchronous, the partially synchronous and the completely asynchronous parts, the third Lyapunov exponent gradually approaches zero and then oscillates around it within a small margin\footnote{with the exception for windows of stability.}. As in the previous case, partial synchronization plays an important role in the appearance of the regime with an additional zero LE.

\begin{figure}[ht]
	\centering{
		\includegraphics[width=\linewidth]{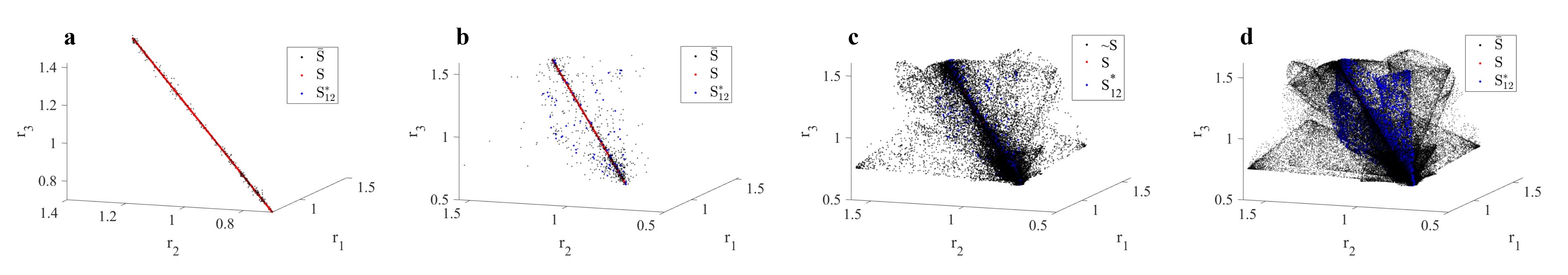}
	}
	\caption{Projections of one of the two components of the Poincar\'e maps of the following attractors on $(r_1,r_2,r_3)$ space: (a) chaotic attractor at $(d/R_0, a/10^6) = (19.7 ,1.4)$; (b) chaotic attractor  $(d/R_0, a/10^6) = (19.75 ,1.4)$; (c) chaotic attractor at $(d/R_0, a/10^6) = (19.86 ,1.4)$; (d) hyperchaotic with and additional zero LE at  $(d/R_0, a/10^6) = (20.15 ,1.4)$.}
	\label{fig:GH_sec}
\end{figure}

The corresponding Poincar\'e maps demonstrating the expansion of the asynchronous and partially synchronous components are presented in Fig.~\ref{fig:GH_sec}\footnote{At each point along the route GH the attractor consists of two large components. Here we show only one of them for better visual clarity.}. In Fig.~\ref{fig:GH_sec}a one can see a small number of asynchronous points located close to $S$. While the number of points in $\bar{S}^*$ and $S^*_{12}$ is still not large in Fig.~\ref{fig:GH_sec}b, a significant amount of these points appear further away from $S$. Fig.~\ref{fig:GH_sec}c shows the chaotic attractor characterized by a substantial presence of the asynchronous component, however the density of the asynchronous points is still the highest close to $S$. The asynchronous part of the attractor visually overshadows $S$ on the projection of the Poincar\'e map at this point. Finally, we show the hyperchaotic attractor with an additional zero LE in Fig.~\ref{fig:GH_sec}d. One can note a clear difference in the density of the points in $S^*_{12}$ and $\bar{S}^*$ that extend far away from $S$.

\begin{figure}[ht]
	\centering{
		\includegraphics[width=\linewidth]{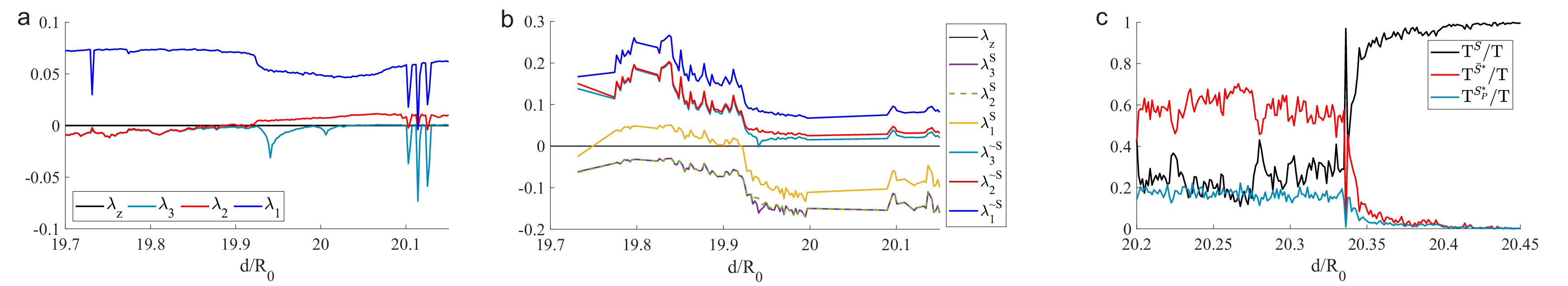}
	}
	\caption{(a) Graph of the largest LEs for the fragment of the path GH: $20.2 < d/R_{0} < 20.45$, (b) the largest LEs, corresponding to the synchronous and asynchronous phases of the the trajectory (when the trajectory spends sufficient time in each phase), (c) the ratios of time the associated with the synchronous, the asynchronous and the partially synchronous components of the attractor.
	}
	\label{fig:GH_split_2}
\end{figure}

The inverse process can be observed afterwards: the gradual increase of the synchronous component of the attractor with simultaneous disappearance of the asynchronous one, accompanying the transition back to synchronous chaotic oscillations ( see Fig.~\ref{fig:GH_split_2}\footnote{the same process, if one considers the decrease of $d/R_0$ starting from $d/R_{0}=20.45$.}).

We think that the same mechanisms is involved in the emergence of the second zero Lyapunov exponent here, as in the path EF. Namely the interplay between completely synchronous, partially synchronous and asynchronous oscillations of the bubbles, creating effects similar to quasiperiodic forcing.

Let us also remark that in some cases (see the discussion above) the existence of strange attractors with an additional zero LEs can be explained by the presence of a codimension-3 bifurcation in the control parameter space \cite{Grines2022, Stankevich2020zero,Kazakov2023}. Generically, for the unfolding a codimension-3 bifurcation we need three independent parameters. In this work we use a fixed value of $\omega$ to reduce the dimensionality of the control parameters space. Varying two parameters ($d/R_0$ and $a$) we do not observe even codimension-2 bifurcations that might indicate the presence of a codimension three bifurcation in the three-dimensional parameter space $(d/R_0, a, \omega)$.

\section{Conclusion}
\label{sec:conlusion}

We have studied oscillations of three encapsulated microbubbles interacting via the Bjerknes force. We have constructed a two-dimensional chart of types of the dynamics in the control parameter space. We have demonstrated that the dynamics of three interacting bubbles may be regular (periodic and quasiperiodic) and chaotic (simply chaotic and hyperchaotic). We have studied several one-parametric routes leading to hyperchaotic oscillations with three positive LEs and chaotic and hyperchaotic oscillations with an additional zero LE. We have analyzed the mechanisms responsible for the emergence of these types of strange attractors. We have discussed two scenarios for the emergence of hyperchaotic attractors with three positive LEs in the considered model. They are based on either the pairs of successive period-doubling bifurcations or Neimark-Sacker bifurcations of unstable periodic orbits (with one-dimensional unstable manifold) lying within a chaotic attractor. Notice also that we have observed the appearance of Lorenz-like quasiattractors \cite{Gonchenko2021, Kaz2024, evgenii2021lorenz} along the routes of the onset of hyperchaotic oscillations. We have found some similarities with the bifurcation mechanisms discussed in \cite{Kazakov2023, Soldatkin2024}. In addition, we believe that the appearance of attractors with an additional zero LE is connected with the existence of partially synchronous oscillations. In this case a sub-cluster of two synchronously oscillating bubbles is influenced by a combination of signals from external driving force and the third bubble. This can create an effect similar to an external quasiperiodic perturbation and, hence, can result in the appearance of an additional zero LE. Finally, let us discuss possible physical and biomedical implications of our results. For observed chaotic attractors with additional zero LE the LEs spectrum has no noticeable oscillations, which might suggest their structural stability. Consequently, it is likely that such attractors can be observed in experiments and the corresponding values of the control parameters can be used in biomedical application, where chaotic dynamics is desirable. Moreover, the prosed bifurcations scenarios do not rely the assumption of completely symmetrical spatial configuration of bubbles. Thus, we believe that these scenarios can be observed in bigger clusters of bubbles.

\section*{Acknowledgments}
This work was supported by the Russian Science Foundation Grant No. 19-71-10048 (Sections 3,4 and 5) and by the Russian Science Foundation Grant No. 19-71-10003 (Section 2). This work were also partially supported by the Basic Research Program at HSE University.

\end{document}